\documentclass{article}

\usepackage{graphicx}
\usepackage[utf8]{inputenc}
\usepackage{amsfonts, amssymb}
\usepackage{amsthm}
\usepackage{amsmath}
\usepackage{enumerate}
\usepackage{graphicx}
\newtheorem{mthm}{Theorem}
\newtheorem{thm}{Theorem}[section]
\newtheorem{mdefn}{Definition}
\newtheorem{defn}{Definition}
\newtheorem{lemma}{Lemma}
\newtheorem{prop}{Proposition}

\newtheorem{cor}{Corollary}
\newtheorem{mcor}{Corollary}
\newtheorem{example}[thm]{Example}
\newcommand{\C}{\mathbb C}

\newcommand{\Om}{\Omega}
\newcommand{\om}{\omega}
\newcommand{\sm}{\setminus}
\newcommand{\ga}{\gamma}

\newcommand{\OO}{{\mathcal O}}
\newcommand{\PP}{{\mathcal{P}}}
\newcommand{\Supp}{\ensuremath{\operatorname{supp}}}

\newcommand{\eps}{\epsilon}
\newcommand{\N}{\mathbb N}
\newcommand{\D}{\mathbb D}
\newcommand{\R}{\mathbb R}

\newcommand{\absv}[1]{\left\lvert{#1}\right\rvert}
\newcommand{\wto}{\stackrel{\mathrm{w}^*}{\to}} % weak star convergence
\DeclareMathOperator{\Cpct}{Cap}
\DeclareMathOperator{\Co}{Co}
\DeclareMathOperator{\Po}{Po}
\newcommand{\ENUM}{\begin{enumerate}}
\newcommand{\ENUMit}{\begin{enumerate}[1.]}
\newcommand{\ENUMA}{\begin{enumerate}[A.]}
\newcommand{\ENUMD}{\begin{enumerate}[D1.]}
\newcommand{\ENUMDp}{\begin{enumerate}[D1'.]}
\newcommand{\ENUMDB}{\begin{enumerate}[DB1.]}
\newcommand{\ENUMDBp}{\begin{enumerate}[DB1'.]}
\newcommand{\ENUMi}{\begin{enumerate}[i)]}
\newcommand{\ENDENUM}{\end{enumerate}}
\newcommand{\ITMZ}{\begin{itemize}}
\newcommand{\ENDITMZ}{\end{itemize}}
\newcommand{\REFEQN}[1]{\begin{equation}\label{#1}}
\newcommand{\ENDEQN}{\end{equation}}
\newcommand{\THM}{\begin{theorem}}
\newcommand{\REFEXA}[1] { \begin{example}\label{#1} }
\newcommand{\ENDEXA}{\end{example}}
\newcommand{\REFMTHM}[1] { \begin{mthm}\label{#1} }
\newcommand{\REFTHM}[1] { \begin{thm}\label{#1} }
\newcommand{\RREFTHM}[2] { \begin{thm}[#1]\label{#2} }
\newcommand{\ENDTHM}{\end{thm}}
\newcommand{\ENDMTHM}{\end{mthm}}
\newcommand{\REFNTH}[1] { \begin{newthm}\label{#1} }
\newcommand{\ENDNTH}{\end{newthm}}
\newcommand{\REFPROP}[1]{\begin{prop}\label{#1} }
\newcommand{\PROP}{\begin{prop}}
\newcommand{\ENDPROP}{\end{prop} }
\newcommand{\DEF}{\begin{defn}}
\newcommand{\REFDEF}[1]{\begin{defn}\label{#1} }
\newcommand{\RREFDEF}[2]{\begin{defn}[#1]\label{#2} }
\newcommand{\ENDDEF}{\end{defn} }
\newcommand{\MDEF}{\begin{mdefn}}
\newcommand{\REFMDEF}[1]{\begin{mdefn}\label{#1} }
\newcommand{\RREFMDEF}[2]{\begin{mdefn}[#1]\label{#2} }
\newcommand{\ENDMDEF}{\end{mdefn} }
\newcommand{\REFLEM}[1]{\begin{lemma}\label{#1} }
\newcommand{\RREFLEM}[2]{\begin{lemma}[#1]\label{#2} }
\newcommand{\LEM}{\begin{lemma}}
\newcommand{\ENDLEM}{\end{lemma} }
\newcommand{\REFCOR}[1]{\begin{cor}\label{#1} }
\newcommand{\COR}{\begin{cor}}
\newcommand{\ENDCOR}{\end{cor} }
\newcommand{\REFMCOR}[1]{\begin{mcor}\label{#1} }
\newcommand{\MCOR}{\begin{mcor}}
\newcommand{\ENDMCOR}{\end{mcor} }
\newcommand{\REFREM}[1]{\begin{Remark}\label{#1} }
\newcommand{\REM}{\begin{Remark}}
\newcommand{\ENDREM}{\end{Remark} }
\newcommand{\REFDEFTHM}[1] { \begin{defthm}\label{#1} }
\newcommand{\ENDDEFTHM}{\end{defthm}}
\newcommand{\eqnref}[1]{(\ref{#1})}
\newcommand{\corref}[1]{Corollary~\ref{#1}}

\newcommand{\defref}[1]{Definition~\ref{#1}}

\newcommand{\lemref}[1]{Lemma~\ref{#1}}
\newcommand{\thmref}[1]{Theorem~\ref{#1}}
\newcommand{\propref}[1]{Proposition~\ref{#1}}

\newcommand{\itemref}[1]{{\emph{\ref{#1}.}}}
\newcommand{\PROOF}{\begin{proof}}
\newcommand{\ENDPROOF}{\end{proof}}

\begin{document}
\title{Convergence of Equilibrium Measures under
$K$-regular Polynomial Sequences and their Derivatives}
\author{
	Henriksen, Christian \\ 
	\texttt{chrh@dtu.dk} \\
	Department of Applied Mathematics and Computer Science \\
	Technical University of Denmark
	\and
	Petersen, Carsten Lunde \\
  \texttt{lunde@math.ku.dk} \\
  Department of Mathematical Sciences\\
  Copenhagen University
	\and
	Uhre, Eva \\
	\texttt{euhre@ruc.dk} \\
	Department of Science and Environment \\
	Roskilde University
}
\date{\today}
\maketitle
\begin{abstract}
  Let $K\subset\C$ be non-polar, compact and polynomially convex.
  We study the limits of equilibrium measures on preimages of compact sets,
  under $K$-regular sequences of polynomials, that center on $K$ 
  and under the sequences of derivatives of all orders of such
  sequences.
  We show that under mild assumptions such limits always exist and
  equal the equilibrium measure on $K$.
  From this we derive convergence of the equilibrium distributions on
  the Julia sets of the sequence of polynomials and their derivatives
  of all orders.
\end{abstract}


\section{Introduction}\label{intro}
There is a longstanding interest in understanding the asymptotic
properties  of particular sequences $(q_k)_k$ of polynomials
$q_k(z) = \ga_k z^{n_k} + \OO(z^{n_k-1})$
of degree $n_k$ and $n_k$ increasing to infinity.
Examples are the classical sequences of orthogonal polynomials, the
Hermite polynomials, the Laguerre polynomials, the Jacobi
polynomials, the Chebyshev polynomials and the Legendre polynomials.
Such polynomials appear in quadrature rules in numerical analysis,
probability theory, representation theory of Lie groups and quantum
groups, and more.
For an arbitrary probability measure $\mu$ on $\C$ with compact
non-polar support $S$, the properties of an extremal sequence
$(q_k)_k= (q_k(z;\mu))_k$
of polynomials that are orthonormal in $L^2(\mu)$,
are studied in the monumental work \cite{StahlTotik} by Stahl and Totik.
Other examples come from holomorphic dynamics.
There one studies the properties of the family of iterates
$(q_k)_k = (P^k)_k$
of a given polynomial $P$ of degree at least $2$ as well as
complex parameter spaces of such polynomials.

A priori, the theory of extremal polynomials and the theory of
holomorphic dynamics seem unrelated.
However pioneering developments have shown asymptotic relations
between the compact support $S$ of a probability measure $\mu$ as
above and the Julia sets and filled-in Julia sets of the sequence of
orthonormal polynomials, see \cite{CHPPOrtho} and \cite{PU}.
In a further development \cite{CHPPCheby} similar relations have been
proven for the $L^\infty(S)$ extremal sequence of dual Chebyshev
polynomials for any non-polar compact set $S\subset \C$.
Furthermore, Bayraktar and Efe generalized such results to so called
asymptotically extremal sequences of polynomials in $L^p(\mu)$, see
\cite{BayraktarEfe}.

In the extremal polynomial setting, a central focus has been to
understand the asymptotic zero distribution of the extremal
polynomials starting with Fej{\'e}r, \cite{Fejer}.
In the holomorphic dynamics setting, the properties of the chaotic
locus, the Julia set, plays a prominent role.
Brolin showed that with minor exceptions, which are well understood,
the preimages under the sequence of iterates of any point in $\C$
equidistribute to the equilibrium measure of the Julia set, see
\cite{Brolin}.
In yet another quite different recent development, Okuyama and Vigny
have shown that with similar exceptions the preimages of any point in
$\C$ under the sequence of derivatives of iterates equidistribute to
the equilibrium measure on the Julia set, see \cite{Okuyama1} and
\cite{OkuyamaVigny}.

In this paper, we begin with the general properties
$K$-centering and $K$-regular (see \defref{C-centering}
and \defref{K-regularity} below),
which are known to hold for both sequences of extremal polynomials
and the sequence of iterates of a fixed polynomial.
We show that such properties are inherited by the sequence of
derivatives, that is, they are passed on to the sequence of
derivatives.
This means for instance that even though a sequence of derivatives of
an extremal sequence of polynomials is not extremal, it still shares
many properties with its parent sequence.
It should be noted here that there is no assumption on the origin of
the sequences satisfying either of the pair of properties $K$-regular
and $K$-centering.
In particular these properties are not tied to any form of
extremality for any probability measure.

This paper is a sequel to the papers \cite{HPU1} and \cite{HPU2} where
we studied inheritance of the asymptotic zero distribution of
sequences of polynomials.
Stahl and Totik, \cite[Theorem 1.1.9.]{StahlTotik},
showed that even in the case of orthonormal polynomials
it is too much to hope for
that the asymptotic zero distribution has a limit.
In this paper, we return to the ideas of \cite{PU} and
\cite{CHPPCheby} and study the convergence properties of balanced
pullbacks of equilibrium distributions on a sequence of compact sets
and the hereditary properties of this convergence.
From this we derive convergence of the equilibrium distributions on
the Julia sets of the sequence of polynomials and their derivatives
of all orders.

\subsection{Definitions and results}
Let $\Omega \subset \C$ be a proper subdomain with non-polar compact
complement  $K:=\C\setminus\Omega$.
Then $K$ is polynomially convex (also known as \emph{full}).
We denote by $J=\partial K=\partial \Omega$ the common topological
boundary, by $g_\Omega$ the Green's function for $\Omega$ with pole
at $\infty$ and by $\omega_K$ the equilibrium measure for $K$.
Then $\Supp(\omega_K) \subseteq J$ with equality if and only if
$\Omega$ is Dirichlet regular.

We denote by $q_k$ a complex polynomial of degree $n_k$:
\begin{equation}\label{meetqk}
q_k(z)=\gamma_{k}z^{n_k}+O\left(z^{n_k-1}\right), n_k \geq 0.
\end{equation}
For a sequence $(q_k)_{k\in \N}$ with $n_k\to\infty$, where $q_k$
is a complex polynomial, we consider two properties.
The first is centering on a compact set $C \subset \C$.
%
% DEFINITION OF CENTERING 
% 
\RREFMDEF{Centering on $C$}{C-centering}
Let $C\subset\C$ be a non-empty and compact set.
We say that a sequence of polynomials 
$(q_k)_{k\in \N}$ \emph{centers} on the set $C$ (or $C$-centers), 
if there exist $N\in \N$ and $R>0$ 
such that the following two conditions hold:
\ENUMit
\item\label{centering_item_1}%[1.] 
  For all  $k\geq N$
  \begin{equation*}%\label{centering1}
    \nonumber q_k^{-1}(\{0\})\subset \D(R).
  \end{equation*}
\item\label{centering_item_2}%[2.] 
  For any closed set $L$ with $L\cap C=\emptyset$
  there exists $M =M(L)>0$ such that 
  \[
    \#(q_k^{-1}(\{0\})\cap L)\leq M\quad \text{for all }k\geq N,
  \]
  where points are counted with their multiplicities as zeroes of $q_k$.
\ENDENUM%\end{itemize}
\ENDMDEF
We remark that we are slightly abusing notation by letting
$\#(q_k^{-1}(\{0\}))$ denote the number of zeroes of $q_k$ counted
with multiplicity and not the cardinality of the set.

The second property is $K$-regularity.
%
% DEFINITION OF K-REGULARITY
%
\RREFMDEF{$K$-regularity}{K-regularity}
  We say that a sequence of polynomials
  $(q_k)_{k\in \N}$ is \emph{$K$-regular} 
  if for some $R>0$
  \begin{equation}\label{genreg}
    \lim_{k\to\infty}{\frac{1}{n_k}}\log \absv{q_k(z)}
    = g_\Om(z)
  \end{equation}
  uniformly on $\C\setminus \D(R)$.
\ENDMDEF

Note that $K$-regularity implies centering condition
\itemref{centering_item_1} with the same $R$.
For the following we shall always assume that $R$ is large enough, so
that $K\subset \D(R)$.

For a polynomial sequence $(q_k)_{k\in \N}$ as in \eqref{meetqk} with
$n_k\to\infty$, we say that a property is \emph{hereditary} if
$(q_k)_k$ having the property implies that the sequence of 
derivatives $(q_k')_{k\in \N}$ also has the property.
In this case we will say that the sequence $(q_k')_{k\in \N}$
\emph{inherits} the property.

Our first result is that the properties $K$-regularity and centering on $K$ are hereditary.
\REFMTHM{Maininheritance}
  Let $K\subset\C$ be non-polar, compact and polynomially convex.
  If a sequence of polynomials $(q_k)_{k\in \N}$ is $K$-regular, then
  so is the sequence of derivatives $(q_k')_{k\in \N}$.
  If a sequence of polynomials $(q_k)_{k\in \N}$ centers on $K$, then
  so does the sequence of derivatives $(q_k')_{k\in \N}$.
\ENDMTHM

Let $q_k^{(m)}$ denote the $m$th derivative of $q_k$. 
It follows from \thmref{Maininheritance} by induction on $m$
that if a sequence $(q_k)_{k\in \N}$ is $K$-regular or centers on $K$, 
then for all $m\geq0$ the sequence
$(q_k^{(m)})_{k\in \N}$ inherits this property. 

\REFMTHM{thmpreimagerDisk}
  Let $K\subset\C$ be non-polar, compact and polynomially convex.
  If a sequence of polynomials $(q_k)_{k\in \N}$ is $K$-regular and
  centers on $K$, then 
  \[
    \forall r>0: \quad \omega_{q_k^{-1}(\overline{\D(r)})}\wto \omega_K.
  \]
\ENDMTHM

We obtain the following as a corollary of Theorems
\ref{Maininheritance} and \ref{thmpreimagerDisk}, since the
properties $K$-regularity and centering on $K$ are hereditary.

\begin{figure}
  \begin{center}
    \begin{tabular}{cc}
      \includegraphics[width=4.8cm]{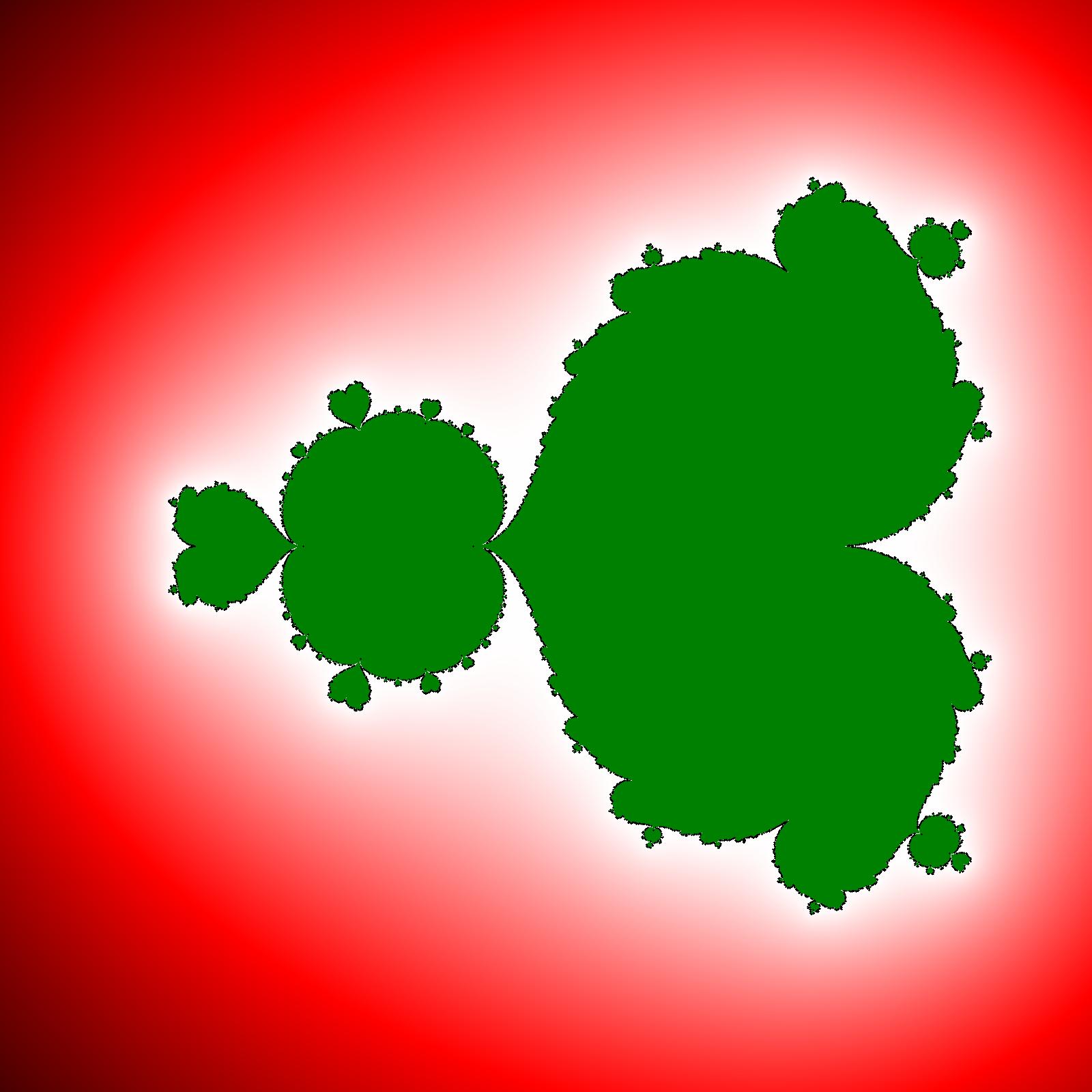} &
    \includegraphics[width=4.8cm]{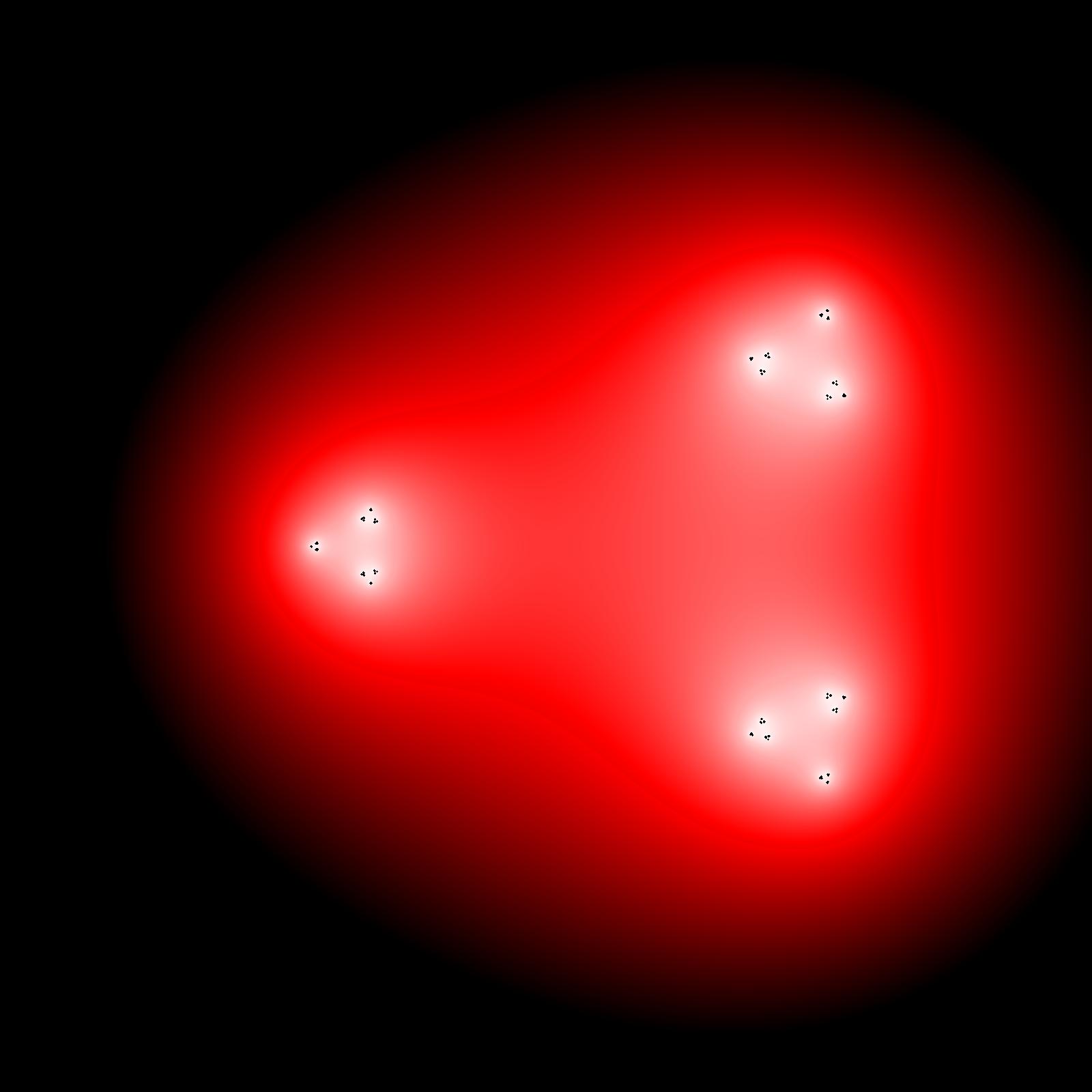}  \\
      \includegraphics[width=4.8cm]{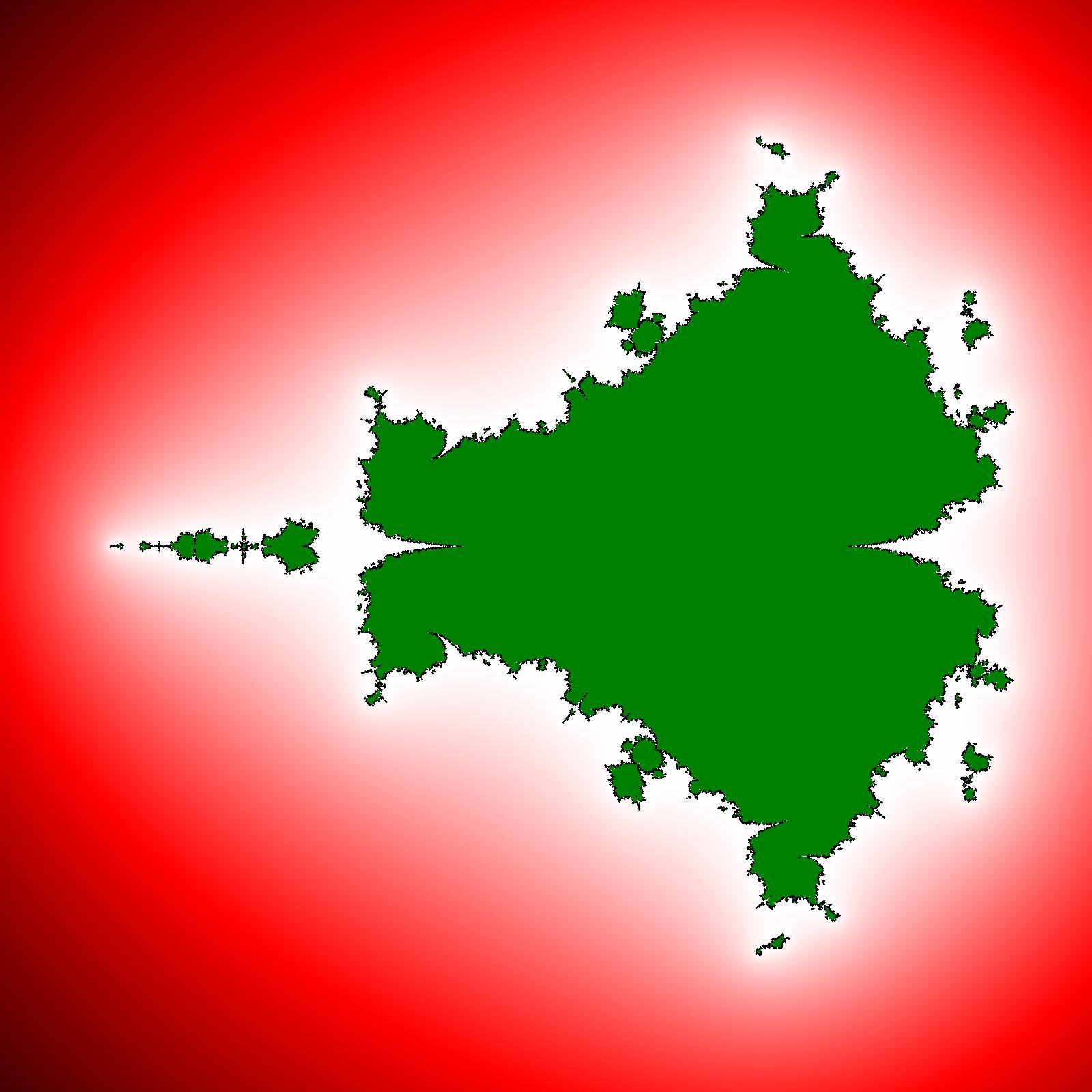} &
    \includegraphics[width=4.8cm]{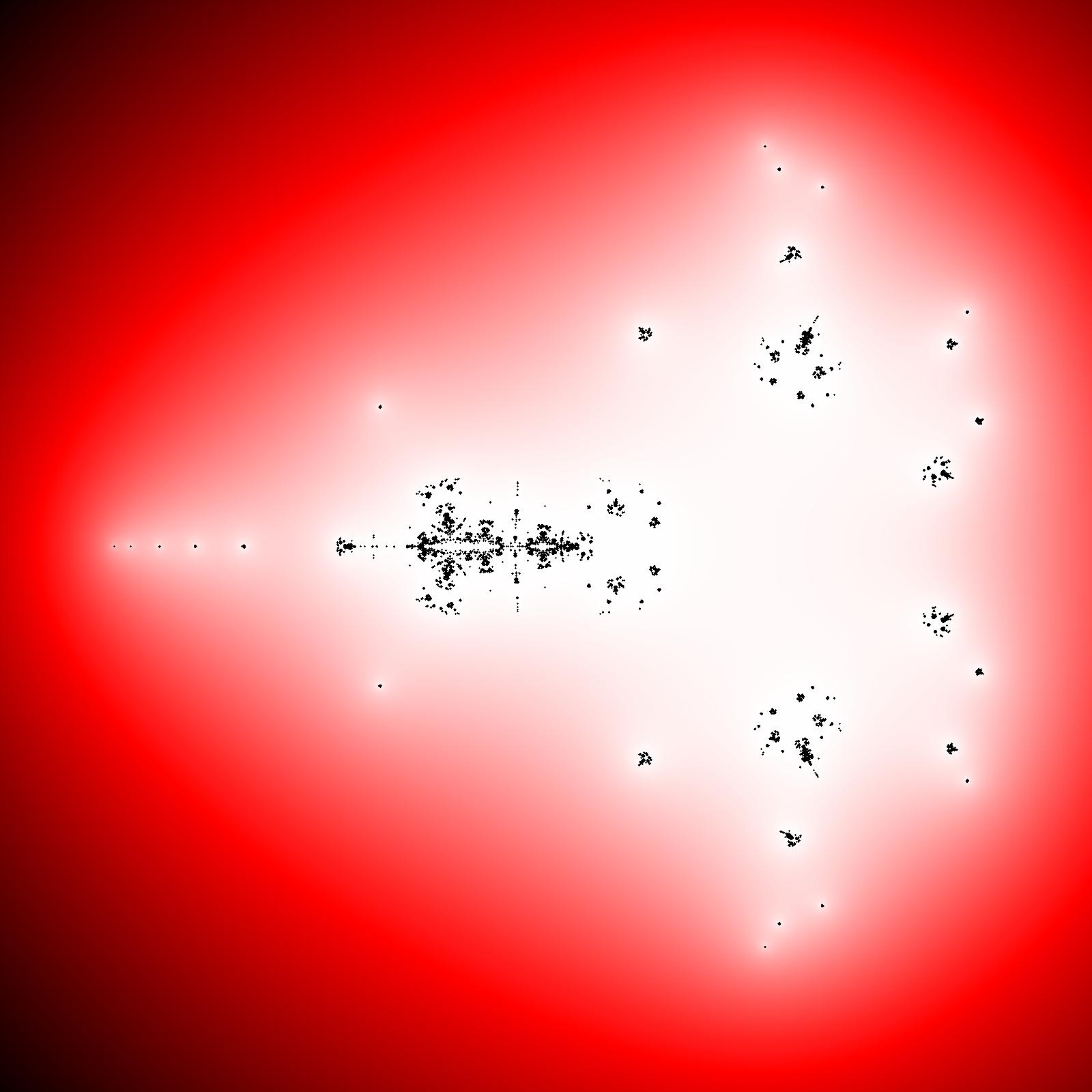}  \\
      \includegraphics[width=4.8cm]{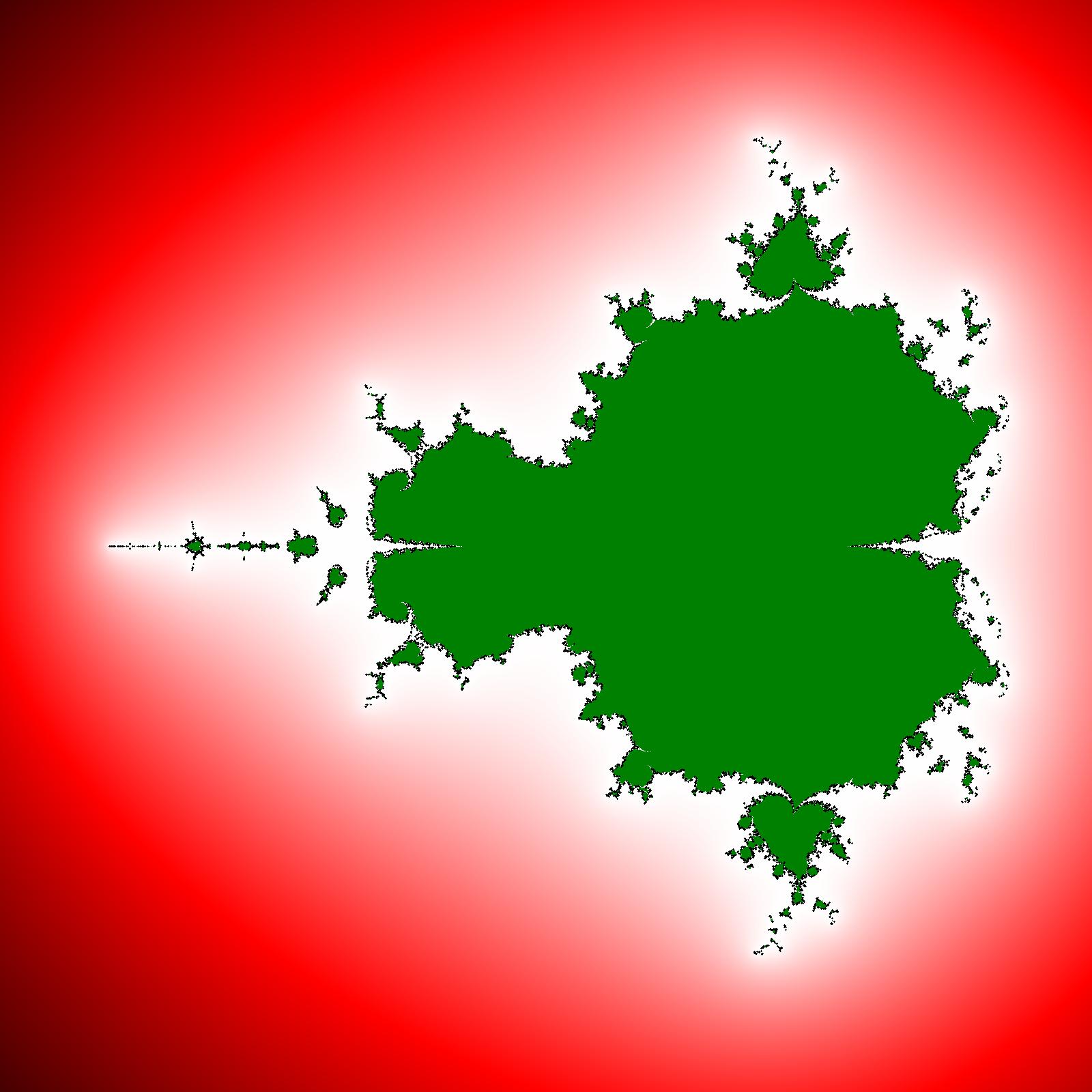} &
    \includegraphics[width=4.8cm]{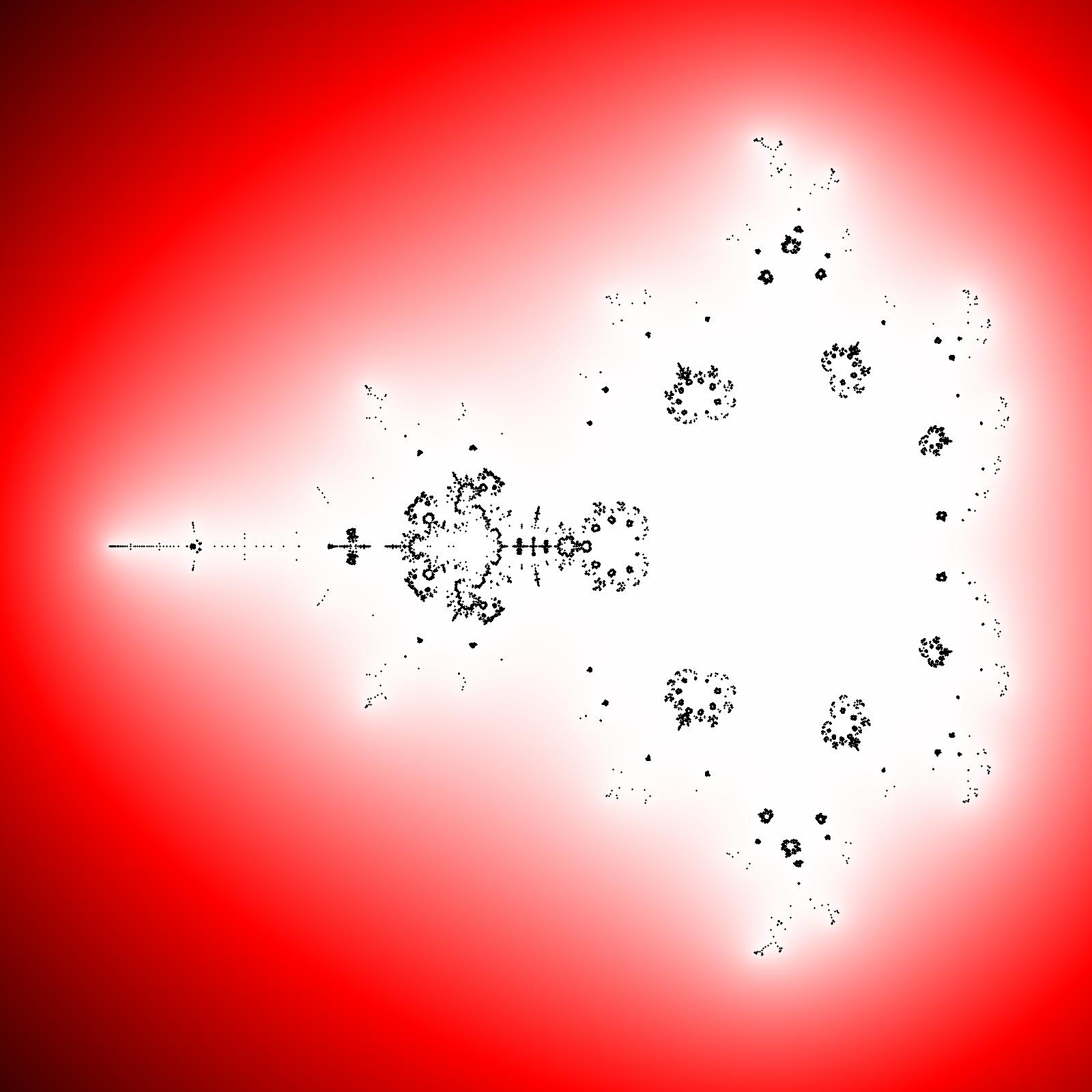}  
    \end{tabular}
  \end{center}
  \caption{ %
  let $P_c(z) = z^2 + c$ and $q_k(c) = P_c^{k}(0)$.
  It follows from \corref{derivatives_qk_equilibrium}
  that $\omega_{K(q_k^{(m)})} \wto \omega_M$, where
  $M$ is the Mandelbrot set.
  Here we illustrate the filled-in Julia set and Green's function
  of $q_k^{(m)}$, for $k = 3, 6, 10$ in the first, second and third row,
  for $m=0$ (left column) and $m=1$ (right column).
  The filled-in Julia set is colored in green, its boundary black,
  amd the Green's function on the complement is illustrated
  by various shades of red.} \label{fig:conv_to_M}
\end{figure}

\MCOR
  Let $K\subset\C$ be non-polar, compact and polynomially convex.
  If a sequence of polynomials $(q_k)_{k\in \N}$ is $K$-regular
  and centers on $K$, then for all $m\geq 1$ 
  \[
    \forall r>0:
    \quad \omega_{(q^{(m)}_k)^{-1}(\overline{\D(r)})}\wto \omega_K.
  \]
\ENDMCOR
For a polynomial $q$ of degree at least $2$, let $q^j$ denote the $j$th iterate of $q$ and 
let $K(q) :=\{z : (q^j(z))_j\textrm{ is bounded}\}$ denote the filled-in Julia set of $q$.
\REFMTHM{TheoremequilibmeasureK}
  Let $K\subset\C$ be non-polar, compact and polynomially convex.
  If a sequence of polynomials $(q_k)_{k\in \N}$ is $K$-regular
  and centers on $K$, then 
  \[
    \omega_{K(q_k)}\wto \omega_K.
  \]
\ENDMTHM
Again we obtain the corresponding result for the sequence of $m$th
derivatives as a corollary of Theorems \ref{Maininheritance} and
\ref{TheoremequilibmeasureK}.
%
% COROLLARY
%
\MCOR\label{derivatives_qk_equilibrium}
  Let $K\subset\C$ be non-polar, compact and polynomially convex.
  If a sequence of polynomials $(q_k)_{k\in \N}$ is $K$-regular and
  centers on $K$, then
  \[
    \omega_{K(q^{(m)}_k)}\wto \omega_K,
  \]
  for $k = 0, 1, \ldots$
\ENDMCOR
To illustrate \thmref{TheoremequilibmeasureK} and
\corref{derivatives_qk_equilibrium}
let $P_c(z) = z^2 + c$ and $q_k(c) = P_c^{k}(0)$.
It is well-known and obvious that $(q_k)_k$ centers on the
Mandelbrot set $M$, and it is shown in \cite{DH82}
that $(q_k)_k$ is $M$-regular.
In particular \corref{derivatives_qk_equilibrium} implies
that $\omega_{K(q_k^{(m)})} \wto \omega_M$.
This is illustrated in Figure \ref{fig:conv_to_M}.

Theorems \ref{thmpreimagerDisk} and \ref{TheoremequilibmeasureK} 
are consequences of the following main technical theorem.
%
% MAIN TECHNICAL THEOREM
%
\REFMTHM{thm:equimeasurecompact}
  Let $K\subset\C$ be non-polar, compact and polynomially convex.
  Let $(L_k)_k$ be a sequence of compact sets in $\C$
  for which 
  $\log(\Cpct(L_k)) = o(n_k)$ as $k \to \infty$ and
  \[
    \exists N\in \N, r_2 > 0\, \forall k \geq N: L_k\subset \D(r_2).
  \]
  If a sequence of polynomials $(q_k)_{k\in \N}$ is $K$-regular
  and centers on $K$, then 
  \[
    \omega_{q_k^{-1}(L_k)}\wto \omega_K.
  \]
\ENDMTHM
\subsection{Examples}
As mentioned above, both centering on $C$ and $K$-regularity are
natural in many contexts.
Here we give more details of some well-known examples.
Many examples will come from sequences of extremal polynomials.
To this end, it is convenient to introduce some notation.
For a closed subset $A\subseteq\C$,  we shall write $\PP(A)$ for the
set of probability measures on $\C$ with compact support in $A$.
Then $\PP(A)$ is sequentially compact under weak* convergence
whenever $A$ is compact. 

\paragraph{Centering.} 
Let $\mu\in \mathcal{P}(\C)$, with non-polar support $S:=\operatorname{supp}(\mu)$,
and let $(P_n)_{n}$ denote the corresponding sequence
of orthonormal polynomials.
Let as usual $\Omega$ denote the unbounded connected component of
$\C\sm S$ and $K=\C\sm\Omega$.
Then the sequence $(P_n)_{n}$ centers on $K$, as was originally shown
by Fejér \cite{Fejer}; see also \cite[Lemma 1.1.3]{StahlTotik} and
\cite[Theorem 1.2]{CHPPOrtho}.

More generally a sequence of extremal polynomials in the sense of Widom
centers on a compact set $K$, see \cite[Lemma 3]{Widom}.

Let $S\subset\C$ be a non-polar compact set with corresponding dual
Chebyshev polynomials $(\mathcal{T}_n)_n$, see \cite{CHPPCheby} for
details.
The by \cite[Section 2.3 and Proposition 4.5]{CHPPCheby}
the sequence $(\mathcal{T}_n)_n$ centers on $K$,
where $K=\Po(S)$ denotes the polynomial convex hull of $S$
(or filled-in $S$). 

Let $S\subset\C$ be a compact set with an infinite number of points;
then an associated sequence of Fekete polynomials $(F_n)_n$ centers on $S$
since in fact all their zeroes belong to $S$,
cf. \cite[Section 5.5]{Ransford}.

Let $q$ denote a polynomial of degree $d>1$ that is not of the form
$\ga z^d$, let $K$ denote its filled-in Julia set, set $\Omega=\C\sm K$,
and let $J=\partial K=\partial \Omega$ denote its Julia set.
Then the sequence of iterates $(q^k)_k$ centers on $J$.
See \cite[Example 6]{HPU1} for details.

\paragraph{$K$-regularity.}
Let $\mu\in \mathcal{P}(\C)$ with non-polar support $S:=\operatorname{supp}(\mu)$
and let $(P_n)_{n}$ denote the corresponding sequence of
orthonormal polynomials.
Let as usual $\Omega$ denote the unbounded connected component of
$\C\sm S$ and let $K=\C\sm\Omega$.
Also, let $\Co(K)$ denote the convex hull of $K$.
Then $K$-regularity of the sequence $(P_n)_{n}$ is equivalent to
$n$th root regularity of the measure $\mu$, that is
$\lim_{n\to\infty}{\frac{1}{n}}\log\absv{P_n(z)} = g_\Om(z)$
locally uniformly on $\C\sm \Co(K)$.

Let $q$ denote a polynomial of degree $d>1$, let $K$ denote its
filled-in Julia set, $\Omega=\C\sm K$, and
$J=\partial K=\partial \Omega$ its Julia set.
Then the sequence of iterates $(q^k)_k$ is $K$-regular, since in fact
$\lim_{k\to\infty}{\frac{1}{d^k}}\log\absv{q^k(z)} = g_\Om(z)$
uniformly on $\Omega$, see e.g.~\cite[Corollary 6.5.4]{Ransford}.

Let $S\subset\C$ be a non-polar compact set with corresponding dual
Chebyshev polynomials $(\mathcal{T}_n)_n$, see \cite{CHPPCheby} for
details.
Let $K=\Po(S)$ denote the polynomial convex hull of $S$
(or filled-in $S$).
Then the sequence $(\mathcal{T}_n)_n$ is $K$-regular, since
$\lim_{n\to\infty}{\frac{1}{n}}\log\absv{\mathcal{T}_n(z)}
= g_\Om(z)$
locally uniformly on $\C\sm \Co(K)$;
see \cite[Theorem 3.9]{Ransford} or \cite[Theorem 3.2]{CSZ1},
whence uniformly on some $\C\setminus\D(R)$ by the uniform
boundedness of zeroes.

Let $S\subset\C$ be a compact set of Capacity $1$, let as usual
$\Omega$ denote the unbounded connected component of $\C\sm S$ and
let $K=\C\sm\Omega$.
Then an associated sequence of Fekete polynomials $(F_n)_n$ is
$K$-regular, see \cite[Theorem 1.18]{Saff}.

Let $(q_k)_k$ be a sequence of asymptotically minimal polynomials 
for some measure $\mu\in\PP(\C)$ as in \cite{BayraktarEfe} and let $K$ 
denote the filled-in support of $\mu$. 
Then $(q_k)_k$ is $K$-regular.

\section{Brief recap of Potential Theory}
Recall that 
the potential $p_\mu$ for $\mu$ a Borel probability measure
is the subharmonic function on $\C$
\REFEQN{potentialdef}
  p_\mu(z) := \int_{\C} \log\absv{z-w}\,d\mu(w) = \log|z| + o(1)
\ENDEQN
and the energy $I(\mu)\in\{-\infty\}\cup\R$ is given by
\[
  I(\mu) := \int_{\C\times\C} \log\absv{z-w}\, d\mu(w)\,d\mu(z)
  = \int_{\C} p_\mu(z)\,d\mu(z).
\]
Here we use the sign convention of Ransford \cite{Ransford}.
Furthermore, for a non-polar compact set $C$,
the capacity $\Cpct(C)$ is given by
\[
  \Cpct(C) = \sup\{\exp(I(\mu)): \mu\in\PP(C)\}
  = \exp(I(\omega_C))
\]
where the unique measure $\omega_C \in\PP(C)$ realizing the
supremum is called
the equilibrium measure.

For $\Omega$ the unbounded connected component of $\C\setminus C$ 
and $K = \C\setminus\Omega$, we have equality of the equilibrium measures 
$\om_C = \om_K = \om_J$ 
where $J$ is the common topological boundary of $K$ and $\Omega$. 

The Green's function $g_\Om := p_{\om_K} - I(\om_K)$ 
is the potential of the equilibrium measure $\omega_K$
normalized to integrate to $0$ 
against $\omega_K$. 

For any $\mu\in\PP(\C)$ and any $R> 0$ such that $\Supp(\mu) \subset \D(R)$ 
it follows immediately from \eqnref{potentialdef} that 
\[
  \forall\; z, \absv{z}\leq R:\; p_\mu(z) \leq \log(2R).
\]
In particular if $q(z) = \ga \prod_{j=1}^n (z-\zeta_j)$ is a polynomial
of degree $n>0$ with $q^{-1}(\{0\}) \subset\D(R)$, then 
\REFEQN{uniform_upper_bound}
  \frac{1}{n}\log\absv{q(z)}
  = \frac{1}{n}\log\absv{\ga} + p_{\mu_q}(z)
  \leq \frac{1}{n}\log\absv{\ga} + \log(2R)
\ENDEQN
on $\overline{\D}(R)$ where
$\mu_q = \frac{1}{n}\sum_{j=1}^{n} \delta_{\zeta_j}$
is the root distribution of $q$.

\section{Inherited properties of the sequence of derivatives}
A property of sequences of polynomials, such as in \eqnref{meetqk} with $n_k\to\infty$, 
is called \emph{hereditary} 
if $(q_k')_{k\in \N}$ has the property whenever $(q_k)_{k\in \N}$ has the property.

Towards a proof of \thmref{Maininheritance}, we define the function
$g_\Omega' : \Omega \to \C$ by letting
$g_\Omega'(z):=2\frac{\partial}{\partial z}g_\Omega(z)
= (\frac{\partial}{\partial x} - i \frac{\partial}{\partial y}) g_\Omega(z)$,
for $z = x+iy\in \Omega$.
Since the Green's function $g_\Omega$ is harmonic in $\Omega$, this is well defined.
In view of the following lemma, $g_\Omega'$ is holomorphic.

\REFLEM{gprime}
Suppose $E \subset \C$ is open and $U: E \to \R$ is harmonic.
Then $2\frac{\partial}{\partial z}U$ is holomorphic.
Furthermore, if $(U_k)_k$ is a sequence of harmonic functions
$U_k: E \to C$ such that $U_k \to U$ locally uniformly,
then $2\frac{\partial}{\partial z}U_k \to 2\frac{\partial}{\partial z}U$
locally uniformly.
\ENDLEM

\begin{proof}
  Since the lemma is a local statement, we can assume $E$ is a disk $\D(z_0, r)$
  and $U_k \to U$ uniformly on $E$.

  Let $V$ be the harmonic conjugate of $U$ normalized such that $V(z_0) = 0$.
  Then $U+iV$ is holomorphic.
  We get
  \[ 
    \frac{d}{dz}(U+iV) = \frac{\partial}{\partial z}(U+iV) = 
    (\frac{\partial}{\partial x} - i\frac{\partial}{\partial y})(U+iV) 
    = 2 \frac{\partial}{\partial z} U,
  \]
  where the last equality is due to the Cauchy Riemann equations. 
  This shows that $ 2 \frac{\partial}{\partial z} U$ is holomorphic.

  Let $V_k$ be the harmonic conjugate of $U_k$ with $V_k(z_0) = 0$,
  and let $f_k = U_k + i V_k - (U + iV)$, $k = 1, 2, \ldots$
  Then $f_k : E \to \C$, $k = 1, 2, \ldots$ form a sequence of holomorphic functions
  with the real part $U_k$ converging uniformly to the zero function
  on $E$.
  If follows by Montel's theorem that $(f_k)_k$ is a normal family.
  Also, any converging subsequence must converge to a function
  $f$ having the real part constantly equal to zero and it follows
  that $f$ is a constant function.
  Since $f_k(0) = 0$ for all $k$, we must have that $f$ is the zero
  function.
  Consequently, the whole sequence $f_k \to 0$ locally uniformly.

  Local uniform convergence of holomorphic functions,
  implies local uniform convergence of the derivatives.
  Hence, $\frac{df_k}{dz} = 2\frac{\partial}{\partial z}(U_k - U) \to 0$
  locally uniformly, concluding the proof of the lemma.
\end{proof}

Since
${\frac{1}{n_k}}\log\absv{q_k(z)}$ is harmonic on
$\C\setminus \D(R)$ for $k$ large enough,
it follows from Lemma \ref{gprime} and the $K$-regularity of
$(q_k)_{k\in \N}$ (equation \eqref{genreg}) that
\[
  \lim_{k\to\infty} {\frac{1}{n_k}}\frac{q_k'(z)}{q_k(z)}
  = g_\Omega'(z)
  =\frac{1}{z}+O\left(\frac{1}{z^2}\right)
\]
locally uniformly on $\C\sm \D(R)$. 

We need the following well-known result. 
\REFPROP{critptsprop}
Let $\mu\in \mathcal{P}(\C)$, with $S:=\operatorname{supp}(\mu)$, then
\begin{equation*}
	\nabla p_\mu\neq 0 \quad \text{on }\C\sm \Co(S),
\end{equation*}
where $\Co(S)$ denotes the convex hull of $S$.
In particular, when $S$ is non-polar and $\Omega$ deontes
the unbounded connected component of $\C\sm S$,
  \begin{equation*}
      g_\Omega'(z)\neq 0	 \quad \text{on }\C\sm \Co(S). 
  \end{equation*}
\ENDPROP

We include a proof for completeness.

\begin{proof}
  We just need to prove $\nabla p_\mu\neq 0$ for every $z$ outside $\Co(S)$,
  since the rest of the proposition then follows easily.
  That $z$ lies outside the convex hull of $S$,
  implies that there exists an open halfplane
  $H$ contaning $S$, such that $z$ is not an element of the closure $\overline{H}$.
  By rotating and translating, we can suppose $H$ is
  the left half plane $\{x + i y: x < 0\}$,
  and $z$ lies on the positive real axis $z = x > 0$.

  Since $S$ is compact, there exists $\epsilon > 0$,
  so that $u < -\epsilon$ for every $w = u + iv \in S$.
  We compute
  $\frac{\partial}{\partial x} \log(\lvert x - w \rvert) 
  = \frac{x-u}{(x-u)^2+v^2} < \frac{1}{\epsilon}$.
  Hence
  $0 < \frac{\partial}{\partial x} \log(\lvert x - w \rvert) < \frac{1}{\epsilon}$
  for every $w \in S$ and $x > 0$.
  In particular, $\frac{\partial}{\partial x}\log(\lvert x - w \rvert)$
  is dominated by the constant function $\frac{1}{\epsilon}$.
  Since $w\mapsto \frac{1}{\epsilon}$ is integrable with respect to $\mu$,
  it follows from the differentiation lemma that
  \[
    \frac{\partial}{\partial x} p_\mu(x) 
    = \frac{\partial}{\partial x} \int_S \log(\lvert x - w \rvert)\, d\mu(w)
    = \int_S \frac{\partial}{\partial x} \log(\lvert x - w \rvert)\, d\mu(w).
  \]
  Since the integrand is strictly positive for all $w \in S$,
  so is $\frac{\partial}{\partial x} p_\mu(x)$.
\end{proof}

\begin{proof}[Proof of \thmref{Maininheritance}]
First, from \propref{critptsprop} we get
$g_\Omega'(z)\neq 0$ on $\C\sm \Co(K)$.
  Combining this with Lemma \ref{gprime} gives
\[
  \frac{1}{n_k}\frac{q_k'(z)}{q_k(z)}
  =g_\Omega'(z)+o(1)
  = g_\Omega'(z)(1+o(1)),
\]
locally uniformly on $\C\sm \D(R)$ as $k\to\infty$.
Moreover,
\begin{equation}\label{q_kprime}
  \absv{q_k'(z)}=n_k\absv{q_k(z)}\absv{g_\Omega'(z)(1+o(1))}	
\end{equation}
so that
\begin{align*}
  \frac{1}{n_k-1}\log\absv{q'_k(z)}
    & = \frac{1}{n_k-1}\left(\log n_k+ \log\absv{q_k(z)}
    +\log \absv{g_\Omega'(z)(1+o(1))}\right) \\ 
    & = \frac{1}{n_k-1}\log\absv{q_k(z)} + o(1) \\
	& = {\frac{1}{n_k}}\log\absv{q_k(z)}(1 + o(1))
\end{align*}
locally uniformly on $\C\sm \D(R)$ as $k\to\infty$. 
Hence,
\begin{equation*} %\label{}
  \lim_{k\to\infty} {\frac{1}{n_k-1}}\log\absv{q'_k(z)}
    =\lim_{k\to\infty} \frac{1}{n_k}\log\absv{q_k(z)},
\end{equation*}
locally uniformly on $\C\sm \D(R)$.

We know from equation \eqref{q_kprime} that that $q_k'$ does not
vanish on $\C\setminus\D(R)$ for $k$ sufficiently large.
Therefore, $(q_k')^{-1}(\{0\})\subset\D(R)$
which means that the difference
$\frac{1}{n_k-1}\log\absv{q'_k(z)}-\frac{1}{n_k}\log\absv{q_k(z)}$
is a (uniformly) bounded harmonic function on $\C\setminus\D(R)$,
which converges to 0 on $\partial\D(R)$.
So, the maximum principle of harmonic functions
implies that the convergence is uniform on $\C\setminus\D(R)$.
This shows that the sequence of derivatives $(q_k')_{k\in \N}$ is
$K$-regular.

That the sequence of derivatives $(q_k')_{k\in \N}$ centers on
$K$ is stated and proved in \cite[Corollary 3]{HPU1}.
\end{proof}

\section{Proofs of convergence of measure}

For the proofs of the theorems we need the following auxilliary results. 
\REFPROP{diagonalsubsequence}
  If a sequence of polynomials $(q_k)_{k\in \N}$ is $K$-regular and
  centers on $K$, then for every subsequence of $(q_k)_{k\in \N}$ there
  exist a further subsequence $(q_{k_j})_{j\in \N}$ and a set $E$
  without accumulation points in $\Omega$ such that
  \REFEQN{extendedreg}
    \lim_{j\to\infty}{\textstyle\frac{1}{n_{k_j}}}\log\absv{q_{k_j}(z)}
    = g_\Om(z)
  \ENDEQN
  locally uniformly on $\Omega\setminus E$.
\ENDPROP
First we need some notation.
For a set $A\subset \C$ we let $d(z,A)=\inf_{a\in A}\absv{z-a}$.
Let
\[
  \Omega_\epsilon
  = \{z\in\Omega :  d(z,K)\geq \epsilon\}.
\]
We introduce the notation $z\succeq w$ to stand for
$d(z,K) \geq d(w,K)$.

\begin{proof}[Proof of \propref{diagonalsubsequence}.]
Let $z_{k,i}$, $i=1, \ldots, n_k$, denote the zeroes of $q_k$ counted
with multiplicity, and numbered with respect to distance to $K$, so
that
\[
  z_{k,1} \succeq z_{k,2} \succeq \cdots \succeq z_{k,n_k}.
\]
We now consider the sequences
$(z_{k,i})_{k\in\N}$ for $i=1, \ldots, n_k$.
Because of centering condition \itemref{centering_item_1}, 
we can assume, passing to a subsequence if necessary,
that $z_{k,1}$ converges to some point $z_1\in\overline{\D}(R)$. 
Possibly passing to a further subsequence, we can obtain
convergence for $z_{k,2}$ simultaneously to a point $z_2$ and so on. 
Passing to the diagonal subsequence, we obtain a subsequence such that 
$z_{k,j}\to z_j$ for each $j\in \N$.

The sequence $(z_j)_{j\in\N}$ of limits respects the initial numbering:
\[
  z_{1} \succeq z_{2} \succeq \cdots
\]

Let $E=\{z_j\}_{j\in\N}$ denote the set of accumulation points of
zeroes obtained in this way.
The set $E$ is without accumulation points in $\Omega$, since for
every $\epsilon>0$ and for every $0<\delta<\epsilon$,
$\#(\Omega_{\epsilon}\cap E)\leq M(\Omega_{\epsilon-\delta})$, where
$M(\Omega_{\epsilon-\delta})$ is the bound on zeroes of $q_k$ in
$\Omega_{\epsilon-\delta}$ from centering condition 2.

We now claim that for the diagonal subsequence constructed above,
\begin{equation}\label{extendedreg2}
  \lim_{j\to\infty}\frac{1}{n_{k_j}}\log\absv{q_{k_j}(z)}
  = g_\Om(z)
\end{equation}
locally uniformly on $\Omega\setminus E$.

To see this, fix $\epsilon>0$ and $\delta>0$ and set
\[
  \Omega_{\epsilon,\delta}=\Omega_{\epsilon}\sm E_\delta,
\]
where $E_\delta=\{z\in\C : d(z,E)<\delta\}$ is an open
$\delta$-neighborhood of $E$.
Choose $N$ such that for $j\geq N$
\[
  q_{k_j}^{-1}(\{0\})\cap \Omega_{\epsilon,\delta/2}
  = \emptyset,
\]
which is possible since the zeroes of $q_{k_j}$ only accumulate
on $E$, by construction of the subsequence.

Let $\Omega^\infty_{\epsilon,\delta}$ denote the unbounded component
of $\Omega_{\epsilon,\delta}$.
Now, $\frac{1}{n_{k_j}}\log\absv{q_{k_j}}$ is harmonic on
$\Omega^\infty_{\epsilon,\delta}$ and uniformly bounded above on
$\D(2R)$ (see \eqref{uniform_upper_bound}),
so the sequence
$(\frac{1}{n_{k_j}}\log\absv{q_{k_j}})_j$ restricted to
$\Omega^\infty_{\epsilon,\delta}\cap \D(2R)$ has a convergent
subsequence, with limit agreeing with $g_\Omega$ on $\D(2R)\sm \D(R)$
by the $K$-regularity property.
Therefore, the limit has to be equal to $g_\Omega$ on
$\Omega^\infty_{\epsilon,\delta}$ by the identity principle for
harmonic functions.
This limit is independent of the subsequence, hence
\[
  \lim_{j\to\infty}{\textstyle\frac{1}{n_{k_j}}}\log\absv{q_{k_j}(z)}
  = g_\Om(z)
\]
uniformly on $\Omega^\infty_{\epsilon,\delta}$ for all
$\epsilon, \delta >0$, which shows \eqref{extendedreg2}.
Compare also with \cite[Lemma 1]{HPU1}.
\end{proof}

We will prove \thmref{thm:equimeasurecompact} first,
and then \thmref{thmpreimagerDisk} and \thmref{TheoremequilibmeasureK}
will follow as corollaries.
		
\begin{proof}[Proof of \thmref{thm:equimeasurecompact}]
We first show that $K$-regularity implies that there exists $N$ so
that  $\omega_{q_k^{-1}(L_k)}\in \mathcal{P}(\overline{\D(R)})$ for
$k\geq N$, whence by sequential compactness of
$\mathcal{P}(\overline{\D(R)})$ with respect to weak* convergence,
any sequence $(\omega_{q_k^{-1}(L_k))})_{k\in\N}$ has a weak*
convergent subsequence with limit in $\mathcal{P}(\overline{\D(R)})$.

By the properties of Green's functions there exists an $\epsilon>0$
so that $g_\Omega(z)\geq \epsilon$ for $\absv{z}\geq R$.
Now choose $N$ so that 
\[
  \absv{\frac{1}{n_{k}}\log\absv{q_{k}(z)} - g_\Om(z)}
  <\epsilon/2 \text{ and }
  \mathrm{e}^{n_k\epsilon/2}>r_2
\]
for $k\geq N$, uniformly for $\absv{z}\geq R$.
Hence 
\[
  \frac{1}{n_{k}}\log\absv{q_{k}(z)}>\epsilon/2
\]
so that 
\[
  \absv{q_k(z)}>\mathrm{e}^{n_k\epsilon/2}>r_2
\]
for $\absv{z}\geq R$ and for $k\geq N$, whence  
\[
  q_k^{-1}(L_k)
  \subset {q_k^{-1}(\overline{\D(r_2)})}\subset \D(R)
  \text{ for } k\geq N,
\]
which means $\omega_{q_k^{-1}(L_k)}\in \mathcal{P}(\overline{\D(R)})$
for $k\geq N$.

Next we show that the limit of every weak* convergent subsequence is
$\omega_K$, which together with the sequential compactness of
$\mathcal{P}(\overline{\D(R)})$ finishes the proof.

Using \cite[Lemma 5.2.5]{Ransford} we have
\[
  \Cpct(q_k^{-1}(L_k))
  =\left (\frac{\Cpct(L_k)}{\absv{\gamma_{k}}}\right)^{\frac{1}{n_{k}}}
\]
and by $K$-regularity of $(q_k)_{k\in \N}$
\begin{equation*}%\label{limitcapacity}
  \lim_{k\to\infty}{\frac{1}{n_{k}}}\log\absv{\gamma_{k}}
  = -\log(\Cpct(K)).
\end{equation*}
Combining the two and using that $\log(\Cpct(L_k)) = o(n_k)$
yields
\begin{equation}\label{therightenergyL_k}
	I(\omega_{q_k^{-1}(L_k)})
    = \log(\Cpct(q_k^{-1}(L_k))
    \to \log(\Cpct(K))
    = I(\omega).
\end{equation}

To see that any limit measure is supported on $K$, let
$\mu \in \mathcal{P}(\overline{\D(R)})$
denote the limit of a weak* convergent
subsequence of $(\omega_{q_k^{-1}(L_k)})_{k\in\N}$.
Then according to \propref{diagonalsubsequence} 
we can pass to a further subsequence, so that  
\begin{equation*}
  \lim_{j\to\infty}\frac{1}{n_{k}}\log\absv{q_{k}(z)} =
  g_\Om(z)
\end{equation*}
locally uniformly on $\Omega\setminus E$, where $E$ is the set
consisting of all the accumulation points of the zeroes of $q_k$
as described in the proof of \propref{diagonalsubsequence}.
For a point $z\in E$ we denote by $m_z$ the number of indices $j$
such that $z_j=z$ and call it the multiplicity of $z$.
Recall that the set $E$ is without accumulation points in $\Omega$,
since $(q_k)_{k\in \N}$ centers on $K$.

Let as before
\[
  \Omega_\epsilon=\{z\in\Omega : d(z,K)\geq \epsilon\}.
\]
Fix $\epsilon>0$ and choose $0<\delta<\epsilon/2$ so that for any
pair of points $z_i\neq z_j \in \Omega_{\epsilon/2}\cap E$,
we have that
$\D(z_i,\delta)\cap \D(z_j,\delta)=\emptyset$.
This means that  $\partial \D(z_j,\delta)\cap E=\emptyset$ for all
$z_j\in \Omega_{\epsilon}$,  and for $k$ large enough all zeroes of
$q_k$ in $\Omega_{\epsilon}$ are inside disks $\D(z_j,\delta)$,
and only those zeroes of $q_k$ converging to $z_j$ are inside the
disk $\D(z_j,\delta)$.
Furthermore, let $\delta$ be chosen small enough
(i.e.
$\delta < 
  \min\{%
    d(z_j,\Omega_{\epsilon})
    : z_j\in E\setminus \Omega_{\epsilon}%
  \}$)
so that
$\D(z_j,\delta)\cap \Omega_{\epsilon}\neq \emptyset$
if and only if $z_j\in \Omega_{\epsilon}$, which is possible because of centering.

Let $z_j\in \Omega_{\epsilon}$.
Since $\partial \D(z_j,\delta)\cap E=\emptyset$ and because of
$K$-regularity (\propref{diagonalsubsequence}), there exists a
number $a>0$ such that
\[
  \frac{1}{n_k}\log\absv{q_k(z)}
  \geq a \quad \text{on } \partial \D(z_j,\delta)
\]
whence
\[
  \absv{q_k(z)} \geq e^{n_k a}>r_2 \quad \text{on } \partial \D(z_j,\delta)
\]
for $k$ sufficiently large.

This means that
\[
  q_k^{-1}(L_k)\cap \Omega_{\epsilon}
  \subset \bigcup_{z_j\in\Omega_{\epsilon}} \D(z_j,\delta).
\]

Fix $z\in E\cap\Om_\eps$ and let $U_{z,k} := q_k^{-1}(\D(r_2))\cap\D(z,\delta)$. 
Then the restriction $q_k : U_{z,k} \to \D(r_2)$ is proper of degree $m_z$, 
where $m_z$ is the multiplicity of $z$ in $E$ as defined above. 
Since $L_k\subset \D(r_2)$ it follows that 
\[
\omega_{q_k^{-1}(L_k)}(\D(z,\delta)) =\frac{m_z}{n_k}. 
\] 

Hence  
\[
\omega_{q_k^{-1}(L_k)}(\Omega_{\epsilon})
\leq \frac{1}{n_k} \sum_{z\in E\cap \Omega_{\epsilon}} m_z
\leq \frac{M(\epsilon/2)}{n_k},
\] 
where $M(\epsilon/2)$ is the asymptotic bound on the number of zeroes of
$q_k$, counted with multiplicity, in $\Omega_{\epsilon/2}$
given by centering on $K$.
Therefore, we have for all $\epsilon>0$
\[
  \omega_{q_k^{-1}(L_k)}(\Omega_{\epsilon}) \to 0
\]
as $k \to\infty$, which by
$\omega_{q_k^{-1}(L_k)}\wto \mu$ 
yields
\[
  \mu(\Omega_{\epsilon})=0 \text{ for all }\epsilon>0,
\]
whence $\operatorname{supp}(\mu)\subseteq K$. 

By \cite[Lemma 3.3.3]{Ransford} and \eqref{therightenergyL_k}:
\[
  I(\mu)
  \geq \limsup_{k\to\infty}I(\omega_{q_k^{-1}(L_k)})
  =I(\omega_K),
\]
so that $\mu=\omega_K$ by uniqueness of the equilibrium measure.
\end{proof}
		
\begin{proof}[Proof of \thmref{thmpreimagerDisk}]
  The theorem is a special case of \thmref{thm:equimeasurecompact},
  with $L_k=\overline{\D}(r)$.
  Clearly $\log(\Cpct(\overline{\D}(r))) = o(n_k)$,
  so the requirement of \thmref{thm:equimeasurecompact} are fulfilled
  for some choice of $r_2 > r$.
\end{proof}

To prove \thmref{TheoremequilibmeasureK}, we need the following lemma.

\begin{lemma}\label{FilledJsetsuniformlybounded}
  If a sequence of polynomials $(q_k)_{k\in \N}$ is $K$-regular and
  centers on $K$, then there exist $R>0$ and $N\in\N$ such that
  \begin{equation*}%\label{Juliasetuniformbound}
    \forall k\geq N:\; K(q_k) \subset q_k^{-1}(\overline{\D(R)})\subset \D(R),
  \end{equation*}
  where $K(q_k)$ is the filled-in Julia set of
  $q_k$.
\end{lemma}

This lemma is analogous to \cite[Proposition 3.3]{CHPPCheby}, and the
proof is essentially the same and left to the reader.
Note that $R$ in this lemma naturally is the same as the $R$ coming from
$K$-regularity and from centering property 1.

\begin{proof}[Proof of \thmref{TheoremequilibmeasureK}]
Again, the theorem is a special case of 
\thmref{thm:equimeasurecompact}, here with
$L_k=K(q_k)$.
By \lemref{FilledJsetsuniformlybounded} there exists $N$ so that 
$K(q_k)\subset \D(R)$ for $k\geq N$. 

Next we have (\cite[Theorem 6.5.1]{Ransford})
\[
  \Cpct(K(q_k))
  =\left (\frac{1}{\absv{\gamma_{k}}}\right)^{\frac{1}{n_{k}-1}}
\]
and by $K$-regularity of $(q_k)_{k\in \N}$
\begin{equation*}%\label{limitcapacityKqk}
  \lim_{k\to\infty}\frac{1}{n_{k}}\log\absv{\gamma_{k}}
  = -\log(\Cpct(K)).
\end{equation*}
Combining the two yields
\begin{equation*}%\label{therightenergyKqk}
  \frac{1}{n_k} \log(\Cpct(K(q_k))) \to 0 \text{ and }
  \Cpct(K(q_k)))\to \Cpct(K).
\end{equation*}
We can thus take $r_2=R$ and fulfill the conditions of
\thmref{thm:equimeasurecompact}.
\end{proof}
%
%
% Declarations section required by Potential Analysis
%
%
\section{Declarations}
\paragraph{Ethical Approval.} Not applicable.
\paragraph{Funding.} All three authors gratefully acknowledge
that they were supported by
\textit{Danish Council for Independent Research~|~Natural Sciences}
via the grant DFF--1026--00267B.
\paragraph{Availability of data and materials.} Not applicable.
\paragraph{Aknowledgement.} The authors thank the referee for useful
comments and suggestions.
%
%
% Bibliography
%
%
\bibliographystyle{plain}
\bibliography{limitmeasure}
\end{document}